\documentclass[journal ]{new-aiaa}
\usepackage[utf8]{inputenc}
\usepackage{textcomp}

\usepackage{graphicx}
\usepackage{amsmath}
\usepackage[version=4]{mhchem}
\usepackage{siunitx}
\usepackage{longtable,tabularx}
\usepackage{todonotes}
\usepackage{tikz, pgfplots}
\usetikzlibrary{positioning, arrows.meta, bending, patterns, shapes}

\title{Generalizing Space Logistics Network Optimization with Integrated Machine Learning and Mathematical Programming}

\author{Koki Ho \footnote{Dutton-Ducoffe Professor, Associate Professor, Daniel Guggenheim School of Aerospace Engineering, AIAA Senior Member} \footnote{Corresponding Author (kokiho@gatech.edu).}, Yuri Shimane \footnote{Ph.D. Candidate, Daniel Guggenheim School of Aerospace Engineering,
Student Member AIAA.}, and Masafumi Isaji \footnote{Ph.D. Candidate, Daniel Guggenheim School of Aerospace Engineering, Student Member AIAA.}}
\affil{Georgia Institute of Technology, Atlanta, GA, 30332}

\begin{document}

\maketitle

\section{Nomenclature}

{\renewcommand\arraystretch{1.0}
\noindent\begin{longtable*}{lll}
$\mathcal{A}$ &\quad=\quad& Set of arcs \\
$\boldsymbol{b}$ &\quad=\quad& Biases of neural networks \\
$\boldsymbol{c}$ &\quad=\quad& Cost coefficient matrix of commodity\\
${c'}$ &\quad=\quad& Cost coefficient matrix of spacecraft\\
$\boldsymbol{d}$ &\quad=\quad& demand vector\\
$g_0$ &\quad=\quad& Gravitational acceleration on Earth \\
$H^\text{payload}$ &\quad=\quad& Concurrency matrix for payload\\
$H^\text{propellant}$ &\quad=\quad& Concurrency matrix for propellant\\
$I_\text{sp}$ &\quad=\quad& Specific impulse \\
$\mathcal{J}$ &\quad=\quad& Objective function \\
$m_d$ &\quad=\quad& Spacecraft structure mass \\
$m_f$ &\quad=\quad& Spacecraft propellant capacity \\
$m_p$ &\quad=\quad& Spacecraft payload capacity \\
$M_\text{UB}$ &\quad=\quad& Upper bound allowed for propellant tank capacity \\
$\mathcal{N}$ &\quad=\quad& Set of nodes \\ 
$Q$ &\quad=\quad& Commodity transformation matrix \\
$\mathcal{T}$ &\quad=\quad& Set of time steps \\
$t_b$ &\quad=\quad& Spacecraft impulsive burn time \\
$\Delta t$ &\quad=\quad& Time of Flight (TOF)\\
$\mathcal{V}$ &\quad=\quad& Set of spacecraft \\
$W$ &\quad=\quad& Launch time window \\ 
$\boldsymbol{W}$ &\quad=\quad& Weights of neural networks \\ 
$\boldsymbol{X}$ &\quad=\quad& Auxiliary variables for neural networks \\
$\boldsymbol{x}$ &\quad=\quad& Commodity flow variable\\
$\boldsymbol{Y}$ &\quad=\quad& Auxiliary continuous variables for neural networks \\
$y$ &\quad=\quad& Spacecraft flow variable \\
$\boldsymbol{Z}$ &\quad=\quad& Auxiliary binary variables for neural networks \\
$\boldsymbol{z}$ &\quad=\quad& Auxiliary binary variables for a quadratic term \\
$\alpha$ &\quad=\quad& Spacecraft structural fraction used in the nonlinear model\\
$\boldsymbol{\epsilon}$ &\quad=\quad& Spacecraft structural coefficient used in the linear model \\
$\boldsymbol{\omega}, \omega$ &\quad=\quad& Auxiliary continuous variables for neural networks \\

\emph{Superscripts/Subscipts}&&\\ 
$i$ &\quad=\quad& Node index (departure) \\
$j$ &\quad=\quad& Node index (arrival) \\
$s$ &\quad=\quad& Neural network layer \\
$t$ &\quad=\quad& Time index \\
$v$ &\quad=\quad& Vehicle index \\
$+$ &\quad=\quad& Outflow from a node \\
$-$ &\quad=\quad& Inflow into a node \\

\emph{Acronyms}&&\\ 
EML &\quad=\quad& Earth-Moon Lagrangian Points \\
GC &\quad=\quad& Gale Crater \\
KSC &\quad=\quad& Kennedy Space Center \\
LEO &\quad=\quad& Low-Earth orbit \\
LLO &\quad=\quad& Low-lunar orbit \\
LMO &\quad=\quad& Low Mars orbit \\ 
LS &\quad=\quad& Lunar surface \\ 
LSP &\quad=\quad& Lunar South Pole \\ 
LTO &\quad=\quad& Lunar transfer orbit \\ 
MILP &\quad=\quad& Mixed-integer linear programming \\
ML &\quad=\quad& Machine learning \\
NN &\quad=\quad& Neural network \\
PAC &\quad=\quad& Pacific Ocean \\
PINN &\quad=\quad& physics-informed neural network \\
PWL &\quad=\quad& Piece-wise linear \\
TOF &\quad=\quad& Time of flight \\
ReLU &\quad=\quad& Rectified linear unit \\
\end{longtable*}}

\section{Introduction}
\lettrine{R}{ecent} growing complexity in space missions has led to an active research field of space logistics and mission design. This research field leverages the key ideas and methods used to handle complex terrestrial logistics to tackle space logistics design problems. A typical goal in space logistics is to optimize the commodity flow to satisfy some mission objectives with the lowest cost. A variety of studies have been conducted to achieve this goal, as reviewed in Ref. \cite{ho2024modeling}.

One of the successful space logistics approaches is network flow modeling and optimization \cite{Ishimatsu2016,Ho2014-Acta,Ho2016}. In this approach, the orbits or surfaces are discretized into nodes, and a network is formed with the arcs connecting these nodes. Each arc is characterized by $\Delta$v and time of flight (TOF), which are pre-computed. A demand is specified at the destination(s) with potentially a time window and the goal is to find a vehicle/commodity transportation routing plan and schedule to satisfy that demand at the lowest cost (e.g., launch mass). Given the linear nature of the network model, the problem is formulated as mixed-integer linear programming (MILP).

A caveat of the conventional MILP-based network approach for space logistics is its incapability of handling nonlinearity. For example, in the MILP formulation, the spacecraft structure mass and fuel/payload capacity are approximated by a linear relationship. However, this oversimplified relationship cannot characterize a realistic spacecraft design \cite{Chen2018,Isaji2022}. 
Other types of nonlinearity can appear when a nonlinear time-dependent trajectory model is considered in an event-driven network, where the time step of each event itself is a variable \cite{Jagannatha2020}.  
The literature that tackled this issue is limited and mostly on an ad hoc basis; for example, Ref. \cite{Chen2018} used a piece-wise linear (PWL) approximation for spacecraft modeling, and Ref. \cite{Jagannatha2020} used a PWL approximation for trajectories. These methods are limited in their applicability by their specific model representation dedicated to their applications. Attempts have also been made to directly solve the mixed-integer nonlinear optimization problem by metaheuristics-based embedded optimization \cite{Taylor2007} or decomposition-based optimization \cite{Isaji2022}, but these approaches either are computationally expensive or require a near-optimal initial guess. These methods also assume that the nonlinear spacecraft or trajectory models are readily available and can be called routinely as part of the optimization process; this is not always the case in practice as we may only have access to a set of data points due to the computational cost of the model or the organizational structure of the design team.

In response to this challenge, this Note develops a new systematic general framework to handle nonlinearity in the MILP-based space logistics formulation using machine learning (ML). Specifically, we replace the nonlinear constraints in the space logistics formulation with trained ML models that are compatible with MILP. The MILP-compatible ML model includes linear regression, PWL approximations, neural networks (NNs) with Rectified Linear Unit (ReLU) activations, decision tree regression, and random forest regression, among others \cite{grimstad2019relu,anderson2020strong,bergman2022janos,biggs2017optimizing}; these models can be translated into MILP formulations by defining additional variables and constraints while maintaining the linearity. This Note provides the first demonstration of using such trained ML models directly in the MILP-based space logistics optimization formulation.\footnote{ML (specifically reinforcement learning) has been used as an outer loop for the MILP-based space logistics optimization \cite{Takubo2022}, but ML models have not been used directly as part of the MILP formulation for space logistics applications.}

The proposed approach has several advantages. First, it generalizes the existing ad-hoc approaches and enables the MILP-based space logistics formulation to handle nonlinearity with a general surrogate model systematically. Second, the resulting solution can provide an effective initial guess to be further refined by gradient-based integrated nonlinear optimization if the nonlinear model is available \cite{Isaji2022}. Third, the proposed method can directly learn from data even when the nonlinear model is not explicitly available. This last advantage can be significant when we need to make decisions based on historical data or computationally expensive models, or when the nonlinear model is managed by a different team and is not readily accessible. We expect this generalized approach to provide a new research direction that leads to a much broader applicability of space logistics network modeling than the state-of-the-art.

\section{Problem Definition: Space Logistics Network Model}
\begin{figure}[h]	
	\centering
	\includegraphics[scale=0.5]{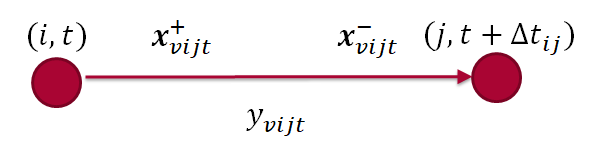}
	\caption{Commodity flow variable definition.}\label{fig:1}
\end{figure}

This section discusses the traditional MILP-based space logistics network model and its challenges. While there are multiple types of formulations for space logistics, we will use the dynamic generalized multi-commodity network flow \cite{Ho2014-Acta} as an example. In this formulation, we consider a time-expanded network represented as graph $\mathcal{G}$ that consists of a set of nodes $\mathcal{N}$ and a set of directed arcs $\mathcal{A}$, over which a set of vehicles $\mathcal{V}$ are moving. Here, the arcs include both the transportation arcs that connect different nodes at different time steps and holdover arcs that connect the same nodes at different time steps. As shown in Fig.~\ref{fig:1}, the variables used to represent the commodity flow over an arc from node $i$ to node $j$ departing at time $t$ for vehicle $v$ is split into outflow $\boldsymbol{x}^+_{vijt}$ and inflow $\boldsymbol{x}^-_{vijt}$. The vehicle flow variable from node $i$ to node $j$ departing at time $t$ for vehicle $v$ is defined as ${y}_{vijt}$. In addition to these flow variables, we also vary the spacecraft parameters including the structure mass $m_{d_v}$, payload capacity $m_{p_v}$, and propellant capacity $m_{f_v}$. The variables and parameters are summarized in Table~\ref{tab_SLvar}.

The resulting optimization problem is given as:
\begin{equation}
\label{SL_obj}
\min \quad \mathcal{J}=\sum_{t\in \mathcal{T}} \sum_{v\in \mathcal{V}} \sum_{(i,j,t)\in \mathcal{A}} ({\boldsymbol{c}_{vijt}}^{T} \boldsymbol{x}^+_{vijt} + {{c'}_{vijt}}^{T} m_{d_v} y_{vijt})
\end{equation}
subject to
\begin{equation}
\label{SL_constr0}
\sum_{v\in \mathcal{V}} \sum_{j:(i,j,t)\in \mathcal{A}} \left[\begin{array}{c}
\boldsymbol{x}^+_{vijt} \\
m_{d_v} y^+_{vijt}
\end{array}\right] - \sum_{v\in \mathcal{V}} \sum_{j:(j,i,t)\in \mathcal{A}} \left[\begin{array}{c}
\boldsymbol{x}^-_{vji(t-\Delta t_{ji})} \\
m_{d_v} y_{vji(t-\Delta t_{ji})}
\end{array}\right]
\leq \boldsymbol{d}_{it} \quad \forall t\in \mathcal{T}\quad \forall i\in \mathcal{N} 
\end{equation}

\begin{equation}
\label{SL_constr1}
\left[\begin{array}{c}
\boldsymbol{x}^-_{vijt} \\
m_{d_v} y_{vijt}
\end{array}\right] = Q_{vijt} \left[\begin{array}{c}
\boldsymbol{x}^+_{vijt} \\
m_{d_v} y_{vijt}
\end{array}\right]\quad \forall t\in \mathcal{T}\quad \forall v\in \mathcal{V} \quad \forall(i, j, t) \in \mathcal{A}
\end{equation}

\begin{equation}
\label{SL_constr2-1}
    {H}^{\text{payload}}_{vij} \boldsymbol{x}^+_{vijt} \leq m_{p_v} y_{vijt}  \quad \forall t\in \mathcal{T} \quad \forall v\in \mathcal{V} \quad \forall(i,j,t)\in \mathcal{A}
\end{equation}

\begin{equation}
\label{SL_constr2-2}
    {H}^{\text{propellant}}_{vij} \boldsymbol{x}^+_{vijt} \leq m_{f_v} y_{vijt}  \quad \forall t\in \mathcal{T} \quad \forall v\in \mathcal{V} \quad \forall(i,j,t)\in \mathcal{A}
\end{equation}

\begin{equation}
\label{SL_vehicle_sizing}
    m_{d_v} = \epsilon_v (m_{d_v} + m_{p_v}) \quad \forall v\in \mathcal{V} 
\end{equation}

\begin{equation}
\label{SL_constr3}
\left\{\begin{array}{ll}
\boldsymbol{x}^\pm_{vijt} \geq \boldsymbol{0}_{p \times 1}  &\text { if } t \in W_{i j} \\
\boldsymbol{x}^\pm_{vijt} =\boldsymbol{0}_{p \times 1} \quad &\text { otherwise }
\end{array} \quad \forall v\in \mathcal{V} \quad \forall( i, j, t) \in \mathcal{A}\right.
\end{equation}

\begin{equation}
\label{xdef}
\begin{aligned}
\boldsymbol{x}_{v i j t}=\left[\begin{array}{c}
x_{1} \\
x_{2} \\
\vdots \\
x_{p}
\end{array}\right]_{v i j t}, \quad \begin{aligned}
&x_{n} \in \mathbb{R}_{\geq 0} \quad \forall n \in \mathcal{C}_{c} \\
&x_{n} \in \mathbb{Z}_{\geq 0} \quad \forall n \in \mathcal{C}_{d}
\end{aligned} \quad \forall v\in \mathcal{V} \quad \forall( i, j, t) \in \mathcal{A} \\
\end{aligned}
\end{equation}

\begin{equation}
\label{udef}
{y}_{v i j t} \in \mathbb{B}  \quad \forall v\in \mathcal{V} \quad \forall( i, j, t) \in \mathcal{A} 
\end{equation}

\begin{equation}
\label{evdef}
\quad m_{p_v}, m_{f_v}, m_{d_v} \in \mathbb{R}_{\geq 0}, \quad \forall v \in \mathcal{V} 
\end{equation}

\begin{table}[h]
\caption{\label{tab_SLvar} Variables and parameters used in the space transportation scheduling problem}
\centering
\begin{tabular}{ll}
\hline\hline Name & 
\begin{tabular}{p{0.7\textwidth}}Description\end{tabular} \\\hline

{\emph{Variables} }&\\

\hline
$\boldsymbol{x}^\pm_{vijt}$ & 
\begin{tabular}{p{0.7\textwidth}}
Outflow (+) and inflow (-) commodity flow variables from node $i$ to $j$ at time $t$ by spacecraft $v$. $\boldsymbol{x}^\pm_{vijt} \geq 0$. Each component of this variable can contain either continuous variables ($\mathcal{C}_{c}$) or discrete variables ($\mathcal{C}_{d}$). This vector will be $p\times1$ vector if the number of commodity types is $p$. 
\end{tabular} \\

$y_{vijt}$ & 
\begin{tabular}{p{0.7\textwidth}}
Spacecraft flow variable, which indicates whether spacecraft $v$ is moving from node $i$ to $j$ at time $t$. This variable is a binary scalar. 
\end{tabular} \\

$m_{p_v}$ & 
\begin{tabular}{p{0.7\textwidth}}
Payload capacity of spacecraft $v$.
\end{tabular}\\

$m_{f_v}$ & 
\begin{tabular}{p{0.7\textwidth}}
Propellant capacity of spacecraft $v$.
\end{tabular}\\

$m_{d_v}$ & 
\begin{tabular}{p{0.7\textwidth}}
Structure mass of spacecraft $v$.
\end{tabular}\\

\hline
{\emph{Parameters}} &\\

\hline
$\boldsymbol{c}_{vijt}$ & 
\begin{tabular}{p{0.7\textwidth}}
Cost coefficient vector for commodities.
\end{tabular}\\

${c'}_{vijt}$ & 
\begin{tabular}{p{0.7\textwidth}}
Cost coefficient for spacecraft. 
\end{tabular}\\

$\boldsymbol{d}_{it}$ & 
\begin{tabular}{p{0.7\textwidth}}
Demands/supplies of different commodities and spacecraft at node $i$ at time $t$. A negative value indicates a demand and a positive value indicates 
a supply.
\end{tabular}\\

$Q_{vijt}$ & 
\begin{tabular}{p{0.7\textwidth}}
Commodity transformation matrix.
\end{tabular}\\

$H^{\text{payload}}_{vij}$, $H^{\text{propellant}}_{vij}$ & 
\begin{tabular}{p{0.7\textwidth}}
Concurrency constraint matrix.
\end{tabular}\\

$W_{ij}$ & 
\begin{tabular}{p{0.7\textwidth}}
Launch window vector, which indicates the available launch window of the spacecraft.
\end{tabular}\\

$\epsilon_v$ & 
\begin{tabular}{p{0.7\textwidth}}
Structural coefficient of spacecraft $v$.
\end{tabular}\\

$\Delta t_{ij}$ & 
\begin{tabular}{p{0.7\textwidth}}
Time of Flight (TOF) from node $i$ to $j$.
\end{tabular}\\

\hline
\emph{Sets} &\\
\hline
$\mathcal{A(N,N,T)}$ & 
\begin{tabular}{p{0.7\textwidth}}
Set of arcs in the time-expanded network.
\end{tabular}\\

$\mathcal{N}$ & 
\begin{tabular}{p{0.7\textwidth}}
Set of nodes.
\end{tabular}\\

$\mathcal{T}$ & 
\begin{tabular}{p{0.7\textwidth}}
Set of time steps.
\end{tabular}\\

$\mathcal{V}$ & 
\begin{tabular}{p{0.7\textwidth}}
Set of spacecraft.
\end{tabular}\\

\hline\hline
\end{tabular}
\end{table}

Equation~\eqref{SL_obj} indicates the objective function, which can be the lifecycle cost or launch mass, depending on the application context. Equation~\eqref{SL_constr0} shows the mass balance constraint between the outflow and inflow at each node.
Equation~\eqref{SL_constr1} shows the mass transformation matrix, where $Q_{vijt}$ is the transformation matrix that indicates the transformation of the commodity over the arc such as propellant consumption.
Equations~\eqref{SL_constr2-1}-\eqref{SL_constr2-2} are the concurrency constraints, enforcing the capacity constraints of the payload and propellant to be below the spacecraft's design parameters $m_{p_v}$ and $m_{f_v}$. 
Equation~\eqref{SL_vehicle_sizing} indicates the structural coefficient constraint to enforce the structural coefficient of a given spacecraft stage. Equations~\eqref{SL_constr3} is the time window constraint, and Eqs. \eqref{xdef}, \eqref{udef}, and \eqref{evdef} show the definitions and domains of each variable.

Note that although Eqs.\eqref{SL_constr2-1}-\eqref{SL_constr2-2} contain quadratic terms $m_{p_v} y_{vijt}$ and $m_{f_v} y_{vijt}$, they can be converted into an equivalent MILP formulation using the big-$\mathcal{M}$ method. For example, $z^{\text{payload}}_{vijt}=m_{p_v} y_{vijt}$ can be converted into the following linear inequalities:
\begin{equation}
    z^{\text{payload}}_{vijt}\leq M y_{vijt}
\end{equation}
\begin{equation}
    z^{\text{payload}}_{vijt}\leq m_{p_v}
\end{equation}
\begin{equation}
    z^{\text{payload}}_{vijt}\geq m_{p_v}-(1-y_{vijt})M
\end{equation}
\begin{equation}
    z^{\text{payload}}_{vijt}\in \mathbb{R}_{\geq 0}
\end{equation}
where $M$ is a large constant such that $m_{p_v} \leq M$.
Recall $y_{vijt}$ is the binary variable representing whether spacecraft $v$ is flying from node $i$ to $j$ at time $t$.
To ensure the payload does not flow without a spacecraft, we introduce the auxiliary variable $z^{\text{payload}}_{vijt}$ such that $z^{\text{payload}}_{vijt} = m_{p_v}$ when $y_{vijt}=1$ (i.e., spacecraft is flying) and $z^{\text{payload}}_{vijt} = 0$ when $y_{vijt}=0$, which are linearly expressed by the above inequalities.
Similarly, $z^{\text{propellant}}_{vijt}=m_{f_v} y_{vijt}$ can also be converted into an equivalent MILP formulation.

While the above formulation is a MILP problem and thus can be rigorously solved, it does not necessarily return a desirable solution given its inherent linear approximation. In fact, any arc constraints from Eqs.~\eqref{SL_constr1}-\eqref{SL_vehicle_sizing} can contain nonlinearity in reality, which is ignored in the above formulation. 
For example, the linear spacecraft model in Eq.~\eqref{SL_vehicle_sizing} based on the structural coefficient is only valid over a limited range of mass and is not a reasonable assumption for global optimization. For a light-weighted spacecraft, one would need a relatively large fraction of structure mass since the mass of necessary subsystems (e.g., avionics, guidance, navigation, control, etc.) does not actually scale down. Thus, in reality, the spacecraft model is often governed by a nonlinear relationship between the structure mass, the payload capacity, and the propellant capacity. Namely, Eq.~\eqref{SL_vehicle_sizing} needs to be replaced with
\begin{equation}
    m_{d_v}=\mathcal{F}(m_{p_v},m_{f_v})
\end{equation}
where $\mathcal{F}$ is a nonlinear function. It is also often the case that this function is not explicitly known and only a set of data points are given (e.g., historical data). In addition to the spacecraft model, the propellant consumption or mass consumption model can also contain a nonlinear relationship or implicit model depending on the application.

The goal of this paper is to develop a general approach to such nonlinearity (including models that are not explicitly known) in the space logistics formulation. The literature on this challenge is limited; one of the few studies that looked into it used a simple customized PWL model to approximate the nonlinear relationship \cite{Chen2018}. In contrast, this paper proposes a more general approach using ML that directly learns the nonlinear model through data and solves the MILP problem with the learned model, expanding the applicability of the network-based optimization model to more general applications.

\section{Theory: Integration of ML Models in Space Logistics Network}
\begin{figure}[h]	
	\centering
	\includegraphics[scale=0.39]{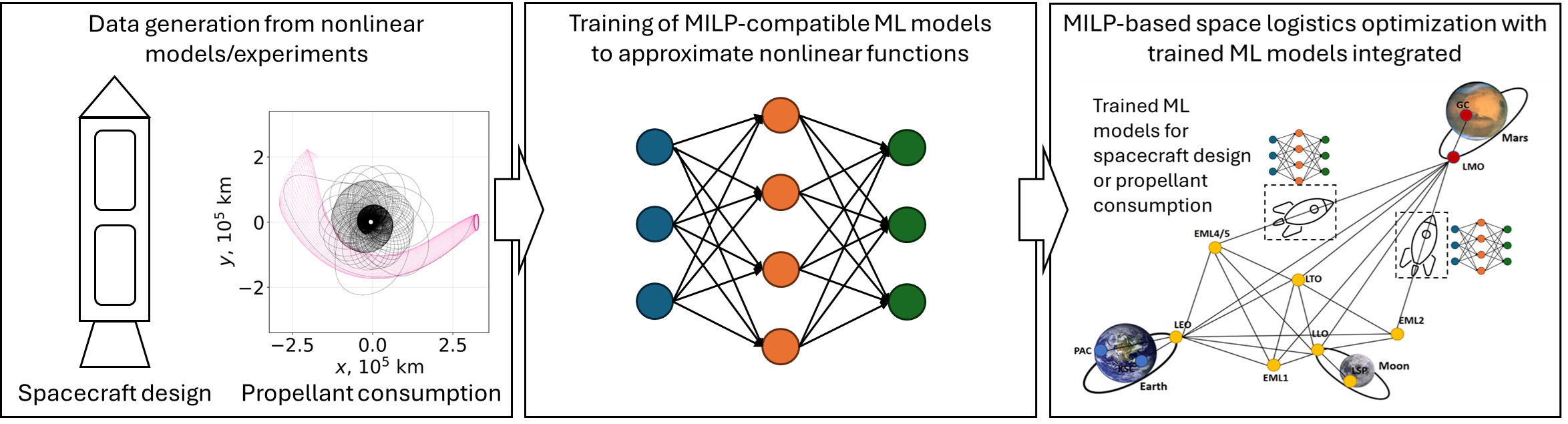}
	\caption{Overall procedure of the proposed method. Part of the figure is adapted from Ref. \cite{Chen2019}.}\label{fig:0}
\end{figure}

This paper proposes leveraging a trained ML model in space logistics network optimization. Specifically, a nonlinear or non-explicit model in the arc constraints Eqs.~\eqref{SL_constr1}-\eqref{SL_vehicle_sizing} can be replaced with an ML model trained beforehand with a set of data points. The chosen ML model needs to be compatible with MILP so it can be written as a set of linear constraints with continuous or integer variables. Such MILP-compatible ML models include linear regression, neural networks with ReLU activations, decision tree, and random forest, among others \cite{grimstad2019relu,anderson2020strong,bergman2022janos,biggs2017optimizing}. Figure \ref{fig:0} shows the overall procedure of the proposed method.

The following shows an example where a complex nonlinear spacecraft model is approximated using a neural network with ReLU activations, which is then converted into linear constraints and integrated into the MILP-based space logistics optimization formulations. The complex spacecraft model represents the relationship between inputs $m_{d_v}$, $m_{p_v}$ and output $m_{f_v}$. We assume that either we have access to that complex nonlinear spacecraft model and can generate a set of training data points $(m_{d_v}, m_{p_v}, m_{f_v})$ from it, or we have access to historical or simulated spacecraft data points $(m_{d_v}, m_{p_v}, m_{f_v})$ that can be used for training. From the training data set, we can train the model:
\begin{equation}
    m_{d_v}=\mathcal{F}_\text{ML}(m_{p_v},m_{f_v})
\end{equation}
For ML, we consider a neural network model that has $S+1$ layers, from $s=0$ to $s=S$, with ReLU activations for the hidden layers and identify activation for the output layer. We assume that each layer $s$ has $n_s$ neurons. The network is fully connected with parameters $\boldsymbol{W}^s$ and $\boldsymbol{b}^s$, where $\boldsymbol{W}^s \in \mathbb{R}^{n_s \times n_{s-1}}$ and $\boldsymbol{b}^s \in \mathbb{R}^{n_{s}}$ for each layer $s$. Thus, the relationship between each layer can be written as follows:
\begin{equation}
    \boldsymbol{X}^s = \sigma{\left(\boldsymbol{W}^s \boldsymbol{X}^{s-1} +  \boldsymbol{b}^s \right)} \quad \forall s=1 \cdots S-1
    \label{relu}
\end{equation}
\begin{equation}
    \boldsymbol{X}^S = \boldsymbol{W}^S \boldsymbol{X}^{S-1} +  \boldsymbol{b}^S
\end{equation}
where $\sigma(\boldsymbol{\omega}):=\max{(\boldsymbol{\omega},0)}$ (component-wise) is the ReLU activation function. The inputs and outputs are
\begin{equation}\label{NN1}
\boldsymbol{X}^0=\begin{bmatrix}m_{p_v}\\m_{f_v}\end{bmatrix}
\end{equation}
\begin{equation}\label{NN2}
\boldsymbol{X}^S=\begin{bmatrix} m_{d_v}\end{bmatrix}
\end{equation}
By using this NN architecture, the resulting MILP formulation involves additional variables $\boldsymbol{X}^s$ for each layer $s$ and additional constraints Eqs.~\eqref{NN1}-\eqref{NN2}; the scalability of the formulation depends on the choice of the NN architecture and its hyperparameters.
 
The ReLU activation function in the form of $Y =\sigma(\omega)=\max{(\omega, 0)}$ can be converted into a set of linear constraints in a variety of ways. A standard and perhaps most intuitive way is as follows \cite{grimstad2019relu}: 
\begin{equation}
Y \geq \omega
\end{equation}
\begin{equation}
Y \leq \omega - M^{-}(1-Z)
\end{equation}
\begin{equation}
Y \leq M^{+}Z
\end{equation}
\begin{equation}
Y \in \mathbb{R}_{\geq 0}
\end{equation}
\begin{equation}
Z \in \mathbb{B}
\end{equation}
where $M^{+}$ and $M^{-}$ are chosen so that $M^{-}\leq \omega \leq M^{+}$.
Note the binary auxiliary variable $Z$ is the binary variable representing whether $\omega$ is greater than 0 or not.
That is, when $Z=1$, the above sets of constraints enforce $Y = \omega$ since $\omega >0$, and similarly, $Y = 0$ is imposed when $Z=0$.
There are recent developments to improve this formulation; see Ref. \cite{anderson2020strong} for an example.

While the above used a neural network for spacecraft design as an example, the approach is general enough to be compatible with a much broader range of problems than the state-of-the-art space logistics formulation. A similar approach can be developed for other MILP-compatible ML models and can also be used to approximate other nonlinearity in the trajectory models or commodity transformation constraints in Eq~\eqref{SL_constr1} as needed.

Several benefits of the proposed approach include the following: First, it can handle general nonlinear models that appear in any constraints in the space logistics formulation. Second, the solution from this method can also be an effective initial guess for gradient-based nonlinear optimization if the original nonlinear model is available. Third, the proposed method can directly learn from a set of training data, and therefore would work even if the original nonlinear model is not readily accessible (e.g., computationally prohibitively expensive or managed by a different team). 

Finally, a few tips on the implementation of ML-based approximations for space logistics are noted. First, due to the stochastic nature of ML model training, it is advisable to test the ML model with an independent set of data and confirm its satisfactory performance before integrating it in MILP. Second, the output from an ML model is not guaranteed to be a positive value, and therefore can potentially lead to an infeasible value. Additional post-processing (such as a ReLU activation) might be needed. Similarly, the output from an ML model can also give a nonzero value even when the inputs are zero. This also justifies the need for post-processing accordingly. Or alternatively, more sophisticated physics-informed NN (PINN) models can be leveraged as well. Finally, the ML model needs to be trained with the dataset in a relevant range of values. When the result requires an ML output that is outside of the range of the training data, the performance would be degraded, and thus re-training would be recommended or an alternative direction needs to be pursued.
    

\section{\label{example}Demonstration: Illustrative Example}
As an illustrative example, we consider a simple problem. We consider a time-expanded network to deliver $1{,}000$ kg of payload mass to the Moon. Our goal is to design the commodity flow and the spacecraft concurrently. This problem is simple enough that its solution can be obtained by solving a system of nonlinear equations. Thus, we can compare the results from the proposed ML-MILP-based approach with the exact solution as a demonstration.

The considered network is shown in Fig.~\ref{LunarMission}, which is extended over a 6-day time horizon with 1-day increments, from $t=0$ to $t=5$ day, and the demand/supply is shown in Table~\ref{LunarMissionDemand}. The objective function is to minimize the initial mass in low-Earth orbit (LEO).

\begin{figure}[hbt!]
\centering
\begin{tikzpicture}[
    node distance=0.25cm and 3cm, 
    arrow/.style={thick, {Latex[length=2mm,width=2mm]}-{Latex[length=2mm,width=2mm]}},
    body/.style={draw, thick, circle, minimum size=1cm, align=center},
    label distance=2mm  
]

  \node[body, fill=blue!30] (PAC) {};
  \node[body, fill=yellow!30, right=of PAC] (LEO) {};
  \node[body, fill=yellow!30, right=of LEO] (LLO) {};
  \node[body, fill=green!30, right=of LLO] (LS) {};

  \node[below=of PAC, align=center] (PAClabel) {Earth};
  \node[below=of LEO, align=center] (LEOlabel) {Low-Earth Orbit\\(LEO)};
  \node[below=of LLO, align=center] (LLOlabel) {Low-Lunar Orbit\\(LLO)};
  \node[below=of LS, align=center] (LSlabel) {Lunar Surface\\(LS)};

  \draw[arrow] (PAC) -- node[above,midway] {TOF= 1 day} (LEO);
  \draw[arrow] (LEO) -- node[above,midway] {TOF= 3 days} node[below,midway] {$\Delta v=4.04$ km/s} (LLO);
  \draw[arrow] (LLO) -- node[above,midway] {TOF= 1 day} node[below,midway] {$\Delta v=1.87$ km/s} (LS);

\end{tikzpicture}
\caption{Lunar campaign network model (modified from \cite{Chen2018}).}
\label{LunarMission}
\end{figure}

\begin{table}[h]
\caption{\label{LunarMissionDemand} Lunar campaign commodity demand and supply}
\centering
\begin{tabular}{cccc}
\hline\hline Commodity Type & Node & Time, day & Supply ($+$) /Demand ($-$)\\\hline
Spacecraft, \# & Earth & 0 & 1 \\
Payload, kg & Earth & 0 & $\infty$ \\
Propellant, kg & Earth & 0 & $\infty$ \\
Payload, kg & Lunar Surface & 5 & -1,000 \\
\hline\hline
\end{tabular}
\end{table}

We consider one spacecraft only for simplicity, and thus the subscript $v$ is dropped. A nonlinear spacecraft model is considered based on the model proposed in Ref. \cite{Taylor2007}, assuming a single-stage LOX/kerosene rocket.
\begin{equation}
    m_d = 2.3931 m_p + \alpha m_f \left(1-\frac{0.2m_f}{M_{UB}}\right)+\frac{0.4189\left(m_f I_\text{sp}g_0/t_b\right)^{0.7764}}{g_0}
    \label{SCDesign}
\end{equation}
where $g_0 = 9.8$ $\text{m/s}^2$ is the gravitational acceleration on Earth; $t_b$ is the spacecraft impulsive burn time that is set to $120$ s; $M_{UB}$ is the upper bound allowed for the propellant tank capacity, which is assumed as $500{,}000$ kg following the original model; and $\alpha$ and $I_\text{sp}$ are the structural fraction and specific impulse, set to $0.045$ and $330$ s, respectively. In this case, since the dependence of $m_d$ on $m_p$ is linear according to the first term of Eq.~\eqref{SCDesign}, we approximate the second and third terms using ML, while leaving the first term as it is:

\begin{equation}
    m_d = 2.3931 m_p + \mathcal{F}_\text{ML}(m_f)
    \label{SCDesign2}
\end{equation}

The model $\mathcal{F}_\text{ML}(m_f)$ is trained as an NN with ReLU activations using a dataset generated by a 1D grid from the range of $m_d\in \left[0,50{,}000\right]$ kg at an increment of $1{,}000$ kg (i.e., $50$ data points). The network is trained using Scikit-learn's MLPRegressor function \cite{scikit-learn} on Python 3.10.9 with 1 hidden layer and 10 neurons and the maximum number of iterations is set to $1{,}000$; the rest of the hyperparameters are set to the default values. With a randomly generated test data set, the $R^2$ value is $0.9907$. With this model, the MILP optimization is performed by Gurobi Machine Learning package \cite{GurobiMachine} on Gurobi 11.0.1.

Table~\ref{results} shows the comparison of the results from the proposed ML-MILP solution with the exact solution, showing that the difference is reasonably small. Due to the stochastic nature of the NN training, the same process is repeated $100$ times with different random seeds; the difference between the proposed solution and the exact solution has a mean of 2.94\% and a median of 0.41\%; the difference between the mean and median is caused by a small number of poorly trained instances, which can be excluded in practice at the testing phase of the NN before being integrated into MILP. The above error is within an acceptable range for preliminary space logistics design. Also, a larger/smaller training data set can potentially result in a more/less accurate result, and a more complex model can represent challenging nonlinearity more accurately at the cost of the computational cost. While the NN-based result is not guaranteed to give a feasible or optimal solution, it can provide a near-optimal initial guess to be fed into the nonlinear model (if it is available) for optimization.  

\begin{table}[h]
\caption{\label{results} Lunar campaign results}
\centering
\begin{tabular}{ccccc}
\hline\hline & Initial Mass in LEO & $m_d$ & $m_p$ & $m_f$ \\\hline
ML(NN)-MILP Solution, kg & 42,941.920 & 5,905.998 & 1,000 & 36,035.923 \\
Exact Solution, kg & 42,811.088 & 5,884.957 & 1,000 & 35,926.131 \\
\hline\hline
\end{tabular}
\end{table}

Note that the purpose of this illustrative example is to show the demonstration of the proposed approach with a simple problem for which we know the exact solution. We do not claim that the chosen particular ML model is the most suitable for this case; rather, our proposal is the general framework that can handle a variety of MILP-compatible ML models.

\section{Conclusion}
This Note provides the first demonstration of using ML to directly incorporate nonlinear models into MILP-based space logistics and mission design optimization. The developed approach has several advantages. First, it can provide a general representation of nonlinearity as a surrogate model in MILP-based space logistics optimization. Second, while the output from the surrogate model itself is not guaranteed to be optimal or feasible for the original nonlinear model, it provides an effective near-optimal initial guess for gradient-based nonlinear optimization if the nonlinear model is available. Third, it is able to learn from data directly even when the nonlinear model is not readily accessible. We believe that this ML-MILP-based framework can significantly generalize the capability of network-based space logistics optimization beyond the state-of-the-art. The specific choice of the ML model would be application-dependent, trading off the training effort, the computational complexity, and the expected performance. We expect the proposed general framework to lead to new future research directions on ML-enhanced space logistics network design for various applications.

\section*{Appendix: Comparison with a Linear Regression Spacecraft Model}
An advantage of the proposed ML-MILP-based approach is that it is general enough to accommodate any type of ML model, not just an NN model. As an example, we can replace the NN model used in Sec.~\ref{example} with a linear regression model. Table~\ref{resultsApp} shows this result. It is shown that the linear regression model also returns a near-optimal solution in this simple example. In general, we expect the NN model to fit better with a nonlinear function but note that neither method is guaranteed to return a feasible or optimal solution; the choice of the method depends on the complexity of the nonlinear function in the original problem.
\begin{table}[h]
\caption{\label{resultsApp} Lunar campaign results: comparison with the linear regression spacecraft mass model}
\centering
\begin{tabular}{ccccc}
\hline\hline & Initial Mass in LEO & $m_d$ & $m_p$ & $m_f$ \\\hline
Linear Regression-MILP Solution, kg & 42,703.819 & 5,867.706 & 1,000 & 35,836.113 \\
Exact Solution, kg & 42,811.088 & 5,884.957 & 1,000 & 35,926.131 \\
\hline\hline
\end{tabular}
\end{table}

\section*{Funding Sources}
This work was conducted with support from the Air Force Office of Scientific Research (AFOSR), as part of the Space University Research Initiative (SURI), under award number FA9550-23-1-0723. Artificial Intelligence (AI) technologies were used in this research for code development and grammar checking.

\bibliography{sample}

\end{document}